\documentclass{elsart}
\input{epsf}

\begin{document}
\parskip 0.05in
\parindent 0.2in

\begin{frontmatter}

\title{Uncertainty relations in models of market microstructure}
\author{Ted Theodosopoulos}
\ead{theo@drexel.edu}
\ead[url]{www.lebow.drexel.edu/theodosopoulos}
\address{Department of Decision Sciences and Department of Mathematics, Drexel University, Philadelphia, PA, USA}


\begin{abstract}
This paper presents a new interacting particle system and uses it as a spin model for financial market microstructure.  The asymptotic analysis of this stochastic process exhibits a lower bound to the contemporaneous measurement of  price and trading volume under the invariant measure in the {\it frozen} phase of the supercritical regime.
\end{abstract}


\begin{keyword}
Spin models, financial market fluctuations, uncertainty relations

\PACS 89.65.Gh
\end{keyword}
\end{frontmatter}
\section{Overview}
This paper presents a stochastic process on a finite periodic lattice.  This process belongs to the family of interacting particle systems \cite{liggett} and is most closely related to the voter process \cite{granovsky}.  However, in contrast to the traditional voter process, the interaction potential for the process presented here incorporates a tunable balance between local and global interactions, characteristic of the minority game and its variants \cite{burgos,challet}.

The proposed stochastic process is interpreted as a spin model for the short-term dynamics of financial markets.  This modeling paradigm is inspired by earlier work of Bornholdt and collaborators \cite{bornholdt1}.  The evolving spin configurations are used to compute a pair of random variables, price and volume, whose variances are investigated.  The goal in this paper is to investigate the behavior of the residual uncertainty in the system, as modeled by the product rms fluctuations of the aggregate price and volume variables.  The approach taken here is comparable to the one Gilmore used in \cite{gilmore}.  The main difference lies with the fact that here the invariant measure of the process is explicitly derived, making unnecessary Gilmore's use of entropy maximization to arrive at the posterior distributions of the observables.  

The proposed stochastic process has two parameters, temperature and a coupling constant.  The process undergoes a natural bifurcation as a function of the coupling constant strength, for decreasing temperatures.  This paper focuses on the analysis of the {\it frozen} phase of the supercritical process.  The invariant measure is characterized and subsequently used to compute the asymptotic variances of the two random variables of interest.  

The impact of perturbations to these asymptotic variances is subsequently investigated.  The stochastic process is extended to include sites whose spins are fixed to a particular value, while the remaining lattice is allowed to relax to its new equilibrium measure.  This perturbation is interpreted as an act of measurement by an exogenous agent who transacts with the market in an attempt to measure simultaneously, with the maximum possible accuracy, the current price and trading volume variables.  This setup is used to formulate the {\it optimal measurement problem}.  The main result of this paper exhibits a strictly positive lower bound to the residual uncertainty in the system for any measurement strategy.

On the one hand, the rigorous analysis of this spin model supplements simulation results based on similar models \cite{bornholdt1,bouchaud4} that generate stylized deviations of real financial markets from the Efficient Markets Hypothesis (EMH) \cite{gabaix1,bouchaud1}.  On the other hand, the framework presented in this paper can be extended to spin models imbedded on more general graphs \cite{majka}. 

The model is described in the next section and the asymptotic dynamics of the unperturbed system are analyzed.  Following that, the measurement problem is presented, leading to the main result.  The paper concludes with a discussion of next steps in this research program.

\section{Modeling Setup and Results}

The state space $X$ of our model is the set of spin configurations on a lattice on the $d$-dimensional torus\footnote{Here we use the notation ${\mathcal T}^d$ to denote the object $\underbrace{ {\mathcal S}^1 \times \ldots \times {\mathcal S}^1}_d$.} $Y \doteq \left({\mathcal Z}/L \right)^d \subset {\mathcal T}^d$, i.e. $X \subset \{-1,1 \}^Y$, for an appropriately chosen $L$ so that $|Y|=N$.  The path of a typical element of $X$ is given by $\eta: Y \times ( 0, \infty ) \longrightarrow \{-1,1\}$ and each site $x \in Y$ is endowed with a (typically $\ell_1$) neighborhood ${\mathcal N} (x) \subset Y$ it inherits from the natural topology on the torus ${\mathcal T}^d$.  In this paper a version of 'rapid stirring' \cite{durrett} is applied by randomizing the neighborhood structure ${\mathcal N}(\cdot)$.  In particular, for each $x \in Y$, ${\mathcal N} (x)$ is a uniformly chosen random subset of $Y$, of cardinality $2d$.  To be more precise, for a set $A$ and a positive integer $k$, let $F(A,k) = \left\{ \left(a_1, a_2, \ldots, a_k \right) \in A^k | a_i \neq a_j \mbox{         } \forall i,j=1,2,\ldots,k \right\}$.  Then, for any $x \in Y$, let $\left\{ {\mathcal N} (x, \cdot) \right\}$ be a family of iid uniform random variables taking values in $F \left( Y \setminus \{x\}, 2d \right)$.

\begin{figure}
\epsfxsize=4in
\epsfbox{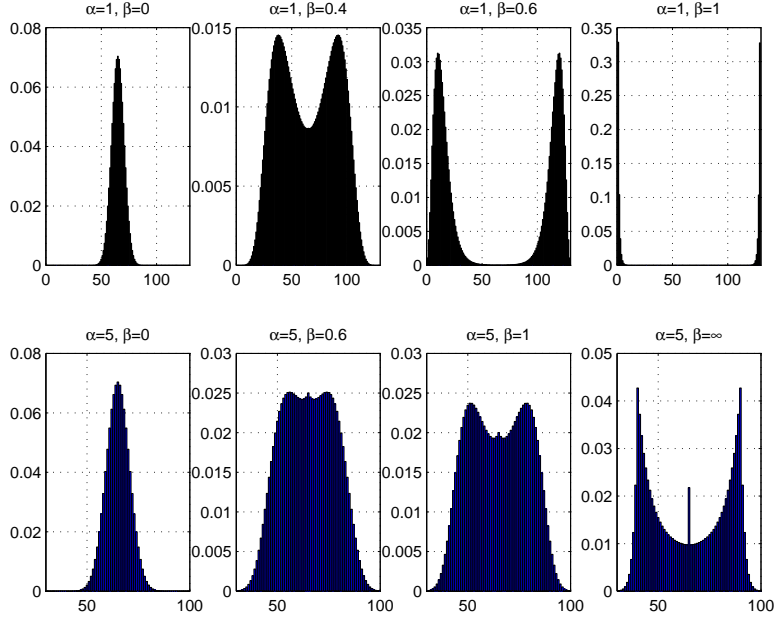}
\caption{The invariant distribution of $N^+$ in the subcritical (top four) and supercritical (bottom four) regime ($d=2, N=128$).}
\label{fig:subsupercritical}
\end{figure}

We construct a continuous time Markov process with transitions occurring at exponentially distributed epochs, $T_n$, with rate 1 \cite{spitzer}.  We proceed to construct a transition matrix for the spins, based on the following interaction potential:
$$h(x,T_n) = \sum_{y \in {\mathcal N}(x,n)} \eta(y,T_n) - \alpha \eta(x,T_n) N^{-1} \left|\sum_{y \in Y} \eta(y,T_n) \right|,$$
where $\alpha>0$ is the coupling constant between local and global interactions.  At time $T_n$ (i.e. the $n$th epoch) a random site $x$ is chosen and its spin is changed to $+1$ with probability $p^+ \doteq \left( 1+ \exp \left\{- 2\beta h \left(x,T_n \right) \right\} \right)^{-1}$ and to $-1$ with probability $p^- = 1- p^+$, where $\beta$ is the normalized inverse temperature.  We define the price and volume processes as follows:
\begin{eqnarray*}
p(t) & = & p^\ast (t) \exp \left\{ \lambda N^{-1} \sum_{y \in Y} \eta(y,t) \right\} \\ 
V(t) & = & N^+ (t) \vee N^- (t)
\end{eqnarray*}
where $a \vee b = \max\{a,b\}$, $a \wedge b = \min\{a,b\}$, $\lambda$ is a parameter, $p^\ast$ is an exogenous previsible `fundamental' price process \cite{bornholdt1} and $N^\pm (t) \doteq \left| \left\{y \in Y \/ | \/ \eta(y,t)=\pm1 \right\} \right|$ denotes the number of sites with a positive or negative spin respectively.
Using the auxiliary variables $X_n = \left| N^+ \left( T_n \right) - N^+ \left( T_{n-1} \right) \right|$ and $Y_n = \left| 2N^+ \left(T_n \right) - N \right|$, we can express $V \left(T_n \right) = \left( Y_n + N \right)/2$ and the volatility, 
$$\sigma \left( \left. \log {\frac {p \left(T_n \right)}{p \left( T_{n-1} \right)}} \right| {\mathcal F}_{n-1} \right) = {\frac {2 \lambda}{N}} \sqrt{{\rm P} \left( X_n = 0 \right) {\rm P} \left( X_n = 1 \right)}.$$
This Markov process is clearly irreducible, and thus possesses a unique invariant measure $\pi_\beta \in {\mathcal M}_1 (X)$ \cite{feller}.  In what follows we consider the $\beta \rightarrow \infty$ limit ({\it frozen} phase) of this process.  The following two observations help characterize the resulting invariant measure:
\begin{itemize}
\item[i.] A site with a spin not in the local majority is unstable and will flip at the next opportunity (i.e. at the next epoch when this site is selected).  When $\alpha \leq 2d$ (what we will refer to as the {\it subcritical regime}) this effect is the only driving force for the process.  Thus, in the long run, every local minority is assimilated and we are led to $\lim_{\beta \rightarrow \infty} \pi_\beta = \left( \delta_0 +\delta_1 \right)/2$, where $\delta_i$ is a Dirac delta mass at $i$ (see the top part of Figure \ref{fig:subsupercritical}).
\item[ii.] When $\alpha > 2d$ (what we will refer to as the {\it supercritical regime}), the earlier effect is supplemented by the fact that a site with a spin in the local majority is unstable (and thus will flip at the next opportunity) for sufficiently high global imbalance (i.e. when $N^+$ is sufficiently high).  Thus, in the supercritical regime (see the bottom part of Figure \ref{fig:subsupercritical}),
\begin{eqnarray*}
\lim_{\beta \rightarrow \infty} \inf \left\{i \in [0,N] | P^\beta_{\scriptscriptstyle --}(i)>0 \right\} & = & N \left( {\frac {1}{2}} - {\frac {d}{\alpha}} \right)  \\
\lim_{\beta \rightarrow \infty} \sup \left\{i \in [0,N] | P^\beta_{\scriptscriptstyle ++}(i)>0 \right\} & = & N \left( {\frac {1}{2}} + {\frac {d}{\alpha}} \right) 
\end{eqnarray*}
where $P^\beta_{\scriptscriptstyle ++}(i)$ is the probability that a $+1$ spin will not flip the next time it is tested when $N^+=i$, at temperature $\beta^{-1}$ ($P^\beta_{\scriptscriptstyle --}(i)$ is the corresponding quantity for a $-1$ spin).
\end{itemize}
\vspace{-0.2in}
\section{The Measurement Problem}

We generalize the process described in the previous section by allowing some randomly chosen sites to be fixed to +1 or -1.  Specifically, let 
\begin{eqnarray*}
& & \pi_\beta \left(i ; k_+, k_- \right) = \lim_{n \rightarrow \infty} {\rm P} \left( N^+ \left(T_n \right) = i \left| \mbox{$k_+$ sites are fixed to +1 and $k_-$ to -1} \right. \right) \\
& & c = \left\lceil \alpha \left| {\frac {i}{N}} - {\frac {1}{2}} \right| \right\rceil , \;\;\;\;\; f_+(i,j) = C_{\scriptscriptstyle j}^{\scriptscriptstyle i-1} C_{\scriptscriptstyle 2d-j}^{\scriptscriptstyle N-i} \left/ C_{\scriptscriptstyle 2d}^{\scriptscriptstyle N-1} \right.  \;\;\;\;\; f_-(i,j) = C_{\scriptscriptstyle j}^{\scriptscriptstyle i} C_{\scriptscriptstyle 2d-j}^{\scriptscriptstyle N-i-1} \left/ C_{\scriptscriptstyle 2d}^{\scriptscriptstyle N-1} \right. \\
& & \begin{array}{cc} P_{\scriptscriptstyle ++}(i) = P_{\scriptscriptstyle --}(i) = 0 & \;\;\;\;\; \mbox{if $i \not\in \left[N \left( {\frac {1}{2}} - {\frac {d}{\alpha}} \right), N \left( {\frac {1}{2}} + {\frac {d}{\alpha}} \right) \right]$} \end{array} \\
& & \begin{array}{cc} P_{\scriptscriptstyle ++}(i) = {\displaystyle \sum_{\scriptscriptstyle j=(d + c) \vee (i+2d-N)}^{\scriptscriptstyle 2d \wedge i}} f_+(i,j) & \;\;\;\;\; \mbox{if $i \in \left[N \left( {\frac {1}{2}} - {\frac {d}{\alpha}} \right), N \left( {\frac {1}{2}} + {\frac {d}{\alpha}} \right) \right]$} \end{array} \\
& & \begin{array}{cc} P_{\scriptscriptstyle --}(i) = {\displaystyle \sum_{\scriptscriptstyle j=0 \vee (i+2d-N)}^{\scriptscriptstyle (d-c) \wedge i}} f_-(i,j) & \;\;\;\;\; \mbox{ if $i \in \left[N \left( {\frac {1}{2}} - {\frac {d}{\alpha}} \right), N \left( {\frac {1}{2}} + {\frac {d}{\alpha}} \right) \right]$} \end{array} \\
& & g(\ell) = C_{\scriptscriptstyle \ell}^{\scriptscriptstyle N-k_+ -k_-} \prod_{\scriptscriptstyle j=0}^{\scriptscriptstyle \ell-1} {\scriptstyle \frac {1- P_{\scriptscriptstyle --} \left(k_+ + j \right)}{1- P_{\scriptscriptstyle ++} \left(k_+ + j +1 \right)}} \;\;\;\;\; \mbox{if $\ell>0$ and $g(0) = 1$} \\
& & Z \left(N, k_+, k_- \right) = \sum_{\scriptscriptstyle i=1}^{\scriptscriptstyle N-k_+-k_-} g(i)
\end{eqnarray*}
where $C_k^n$ denotes the combinations $n$ choose $k$.

\newtheorem{propo}{Proposition}
\begin{propo}
For any choice of $0 \leq k_+ + k_- < N$, the invariant measure of the generalized process is given by:
$$\lim_{\beta \rightarrow \infty} \pi_\beta \left( \ell ; k_+, k_- \right) = \left\{ \begin{array}{ll} {\frac {g(\ell)}{ 1+ Z \left(N, k_+, k_-\right)}} &
\mbox{for $\ell=0, 1, \ldots, N- k_+ -k_-$} \\
0 & \mbox{otherwise} \end{array} \right.$$
\label{propo:invmeas}
\end{propo}

\begin{figure}
\epsfxsize=4in
\epsfbox{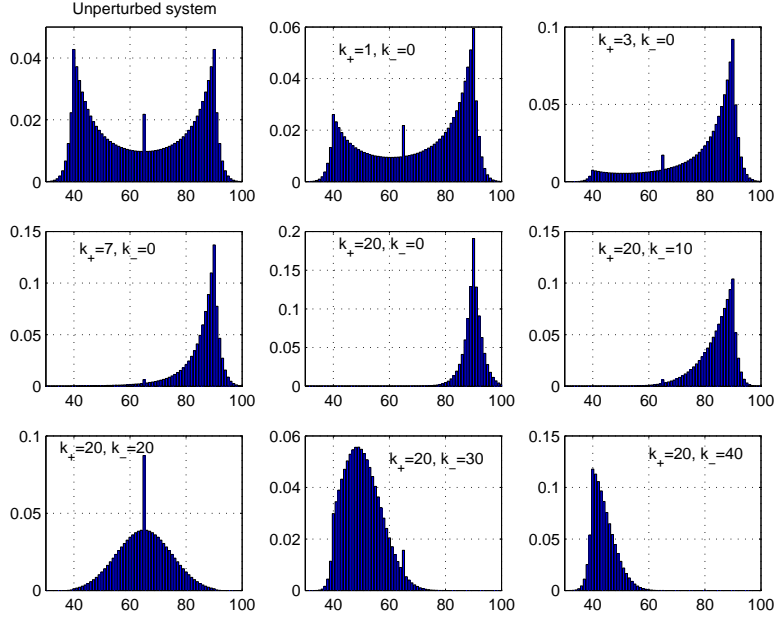}
\caption{This graph shows the invariant distribution for a supercritical process in the {\it frozen} phase with a number of fixed spins.  In all cases $\alpha=5$, $d=2$ and $N=128$.}
\label{fig:kplkmi}
\end{figure}

\begin{figure}
\epsfxsize=4in
\epsfbox{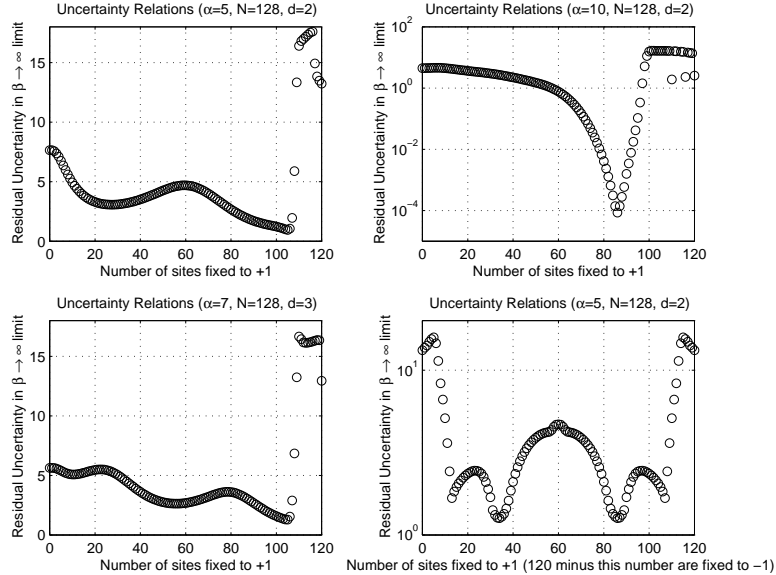}
\caption{This graph shows the effect of fixing the spins for a subset of the sites to the residual variance for a supercritical process in the {\it frozen} phase.  The x-axis for the top two and bottom left figures represents the number of randomly chosen sites (out of a total of $N=128$) that were fixed to +1.  Finally, the bottom right figure shows the result of always fixing 120 randomly chosen sites (out of the total 128), with the x-axis representing the number of those sites that are fixed to +1 (in each case the remaining ones are fixed to -1).}
\label{fig:measurement}
\end{figure}

Figure \ref{fig:kplkmi} shows how the invariant measure of the process becomes skewed when $k_+ \neq k_-$.  Let $\pi_\infty \doteq \lim_{\beta \rightarrow \infty} \pi_\beta$.  Using the measure described in Proposition \ref{propo:invmeas} we obtain the following asymptotic expressions for the two auxiliary random variables $X_n$ and $Y_n$:

\begin{propo}
In the {\it frozen} phase of an unperturbed ($k_+ = k_- = 0$) supercritical process, the following asymptotic moment expressions hold:
\begin{eqnarray*}
& & \lim_{\beta \rightarrow \infty} \lim_{n \rightarrow \infty} {\bf E}^\beta \left[ X_n \right] = 0   \\ 
& & \lim_{\beta \rightarrow \infty} \lim_{n \rightarrow \infty} {\bf E}^\beta \left[ X_n^2 \right] = 1 - {\frac {2}{N}} \sum_{i=0}^N i \pi_\infty \left( i ; 0, 0 \right) P_{++} (i) \\
& & \lim_{\beta \rightarrow \infty} \lim_{n \rightarrow \infty} {\bf E}^\beta \left[ Y_n \right] = N \left(1 - \pi_\infty \left( N/2; 0, 0 \right) \right) - 4 \sum_{i=0}^{{\frac {N}{2}}-1} i \pi_\infty \left( i; 0, 0 \right) \\
& & \lim_{\beta \rightarrow \infty} \lim_{n \rightarrow \infty} {\bf E}^\beta \left[ Y_n^2 \right] = 4 \sum_{i=0}^N i^2 \pi_\infty \left( i; 0, 0 \right) - N^2 
\end{eqnarray*}
\label{propo:moments}
\end{propo}
Finally, Figure \ref{fig:measurement} shows the residual uncertainty, as measured by the product of the asymptotic standard deviations of $X_n$ and $Y_n$ in the {\it frozen} phase, as we perturb $k_+$ and $k_-$ away from 0.  Observe that in all cases the residual uncertainty is lower bounded strictly away from 0, and the (possibly multiple) minima of the residual uncertainty are attained for intermediate-sized perturbations.

\section{Conclusions and Next Steps}

This short paper describes a new interacting particle system and investigates its asymptotic behavior.  It is shown that in the supercritical regime the {\it frozen} phase includes complex evolving patterns of coexistence of both spins.  This setup is used to formulate the problem of measuring a pair of aggregate random variables defined on the spin configurations and it is shown that the minimum attainable residual measurement uncertainty is strictly positive.  

This stochastic process is interpreted as a spin model for market microstructure.  In this regard, it is related to the model proposed earlier by Bornholdt \cite{bornholdt1}, and variants of the voter process and the minority game \cite{burgos,challet}.  As a next step in this research program, the difference between these models and the one proposed here needs to be evaluated.  Furthermore, the extension of the asymptotic study presented here to the finite temperature case and the estimation of the convergence rate to the invariant measure presented in Proposition \ref{propo:invmeas} should be pursued.  Subsequently, the process presented here can be generalized by being imbedded in more general graphs, and a {\it thermodynamic limit} can be formulated, by allowing the number of agents to approach infinity in an appropriate manner.  

\section{Acknowledgements}
The author would like to thank Ming Yuen and an anonymous referee for their careful reading and insightful comments which led to a cleaner presentation.

\end{document}